\documentclass[11pt]{article} 
\usepackage{amsmath,amssymb}

\input{prepictex.tex} 
\input{pictex.tex}
\input{postpictex.tex}

\def\C{\mathbb C}
\def\R{\mathbb R}
\def\T{\mathbb T}
\def\Z{\mathbb Z}

\def\Aut{\operatorname{Aut}}
\def\coker{\operatorname{coker}}
\def\Pol{\operatorname{Pol}}

\newtheorem{prop}{Proposition}
\newtheorem{coro}[prop]{Corollary}

\begin{document}

\title{Quantum Lens Spaces and Principal Actions \\ on Graph $C^*$-Algebras}

\author{Wojciech Szyma\'nski\\ \\
School of Mathematical and Physical Sciences\\
The University of Newcastle, Callaghan, NSW 2308, Australia\\
e-mail: wojciech@maths.newcastle.edu.au}

\date{24 September, 2002}

\maketitle

\abstract{We study certain principal actions on noncommutative 
$C^*$-algebras. Our main examples are the $\Z_p$- and $\T$-actions on the 
odd-dimensional quantum spheres, yielding as fixed-point algebras 
quantum lens spaces and quantum complex projective spaces, respectively. 
The key tool in our analysis is the relation of the ambient $C^*$-algebras 
with the Cuntz-Krieger algebras of directed graphs. A general 
result about the principality of the gauge action on 
graph algebras is given.} 

\addtocounter{section}{-1}

\section{Introduction} 

Classical complex projective spaces may be defined as orbit 
spaces of free actions of the circle group on odd-dimensional spheres. 
This idea can be extended to the noncommutative world. Namely,  
the $C^*$-algebras $C(\C P^{n-1}_q)$ of the 
quantum complex projective spaces are defined as fixed-point algebras 
for certain $\T$-actions on the quantum odd-dimensional spheres 
$C(S^{2n-1}_q)$ \cite{vs}. (Note that according to a recent result of  
Hawkins and Landi \cite{hl}, the Vaksman-Soibelman quantum 
spheres coincide with the ones defined by Reshetikhin, Takhtadzhyan 
and Faddeev \cite{rtf}.) Similarly, 
classical (generalized) lens spaces may be defined as orbit spaces 
of free actions of finite cyclic groups on odd-dimensional spheres. 
This construction was extended to noncommutative setting in \cite{hs2}. 
Therein, natural $\Z_p$-actions on the odd-dimensional 
quantum spheres are considered whose 
fixed-point algebras constitute the $C^*$-algebras of continuous 
functions on the quantum lens spaces. This definition gives rise to 
algebras which are in general non-isomorphic with the ones 
constructed earlier by Matsumoto and Tomiyama \cite{mt}. 

In this note, we show that the above-mentioned actions on 
noncommutative $C^*$-algebras are principal in the sense of Ellwood 
\cite{e}. Special cases include the $\T$-action on $C(SU_q(2))$ 
yielding the standard Podle\'s sphere $C(S^2_{q0})$ \cite{p}, 
and the $\Z_2$-actions on $C(S^{2n-1}_q)$ giving rise to the 
quantum odd-dimensional real projective spaces 
$C(\R P^{2n-1}_q)$ \cite{hs1}. 

In order to smoothly handle $C^*$-algebraic complications we 
make use of the machinery of Cuntz-Krieger 
algebras of directed graphs. In fact, even though the 
$C^*$-algebras under consideration are defined via complicated relations, 
all of them are isomorphic to certain graph $C^*$-algebras \cite{hs1,hs2}. 
This makes the determination of their ideal structure and $K$-theoretic 
invariants a matter of routine calculations (cf. \cite{hs1,hs2}). It also 
greatly simplifies the proofs of principality of the actions. Indeed, both 
the $\Z_p$-actions giving the $q$-lens spaces and the $\T$-actions giving 
the quantum complex projective spaces are defined via very simple 
formulae when transported to the graph algebras. In fact, the latter 
become nothing else but the 
canonical gauge actions of the corresponding graph algebras. This last 
observation motivates our brief discussion of principality of the 
gauge actions on arbitrary graph algebras.

\section{Graph $C^*$-algebras} 

We briefly recall the concept of a graph $C^*$-algebra. (For more details we refer 
the reader to \cite{kprr} and \cite{bprs}.) Let $E$ be a countable graph with 
the set of vertices $E^0$ and the set of directed edges $E^1$. (If $e\in E^1$ 
then $s(e)$ is the source of $e$ and $r(e)$ is its range.) For simplicity sake we 
assume that every vertex in $E$ emits only finitely many edges. Then $C^*(E)$ is 
defined as the universal $C^*$-algebra generated by partial isometries 
$\{S_e\;|\;e\in E^1\}$ with mutually orthogonal ranges and by projections 
$\{P_v\;|\;v\in E^0\}$ such that 
\begin{eqnarray} 
S_e^* S_e & = & P_{r(e)} \mbox{ for each } e\in E^1, \\ 
P_v & = & \sum_{s(f)=v}S_f S_f^* \mbox{ for each } v\in E^0 
\mbox{ emitting at least one edge.} \label{vertex} 
\end{eqnarray} 
These graph algebras generalize and contain as a subclass the 
classical Cuntz-Krieger algebras \cite{ck,c}. 

According to the general Cuntz formula valid for graph algebras 
(cf. \cite[Theorem 3.2]{rs}), the $K$-theory of $C^*(E)$ can be 
calculated as follows. Let $E^0_+$ be the set of those vertices of 
$E$ that emit at least one edge, and let $\Z E^0$ and $\Z E^0_+$ be the 
free abelian groups with generators $E^0$ and $E^0_+$, respectively. Let 
$A_E:\Z E^0_+\longrightarrow\Z E^0$ be the map defined by 
\begin{equation}
A_E(v)=\left(\sum_{s(e)=v}r(e)\right)-v.
\end{equation} 
Then 
\begin{equation}\label{k}
K_0(C^*(E))\cong\coker(A_E), \;\;\; K_1(C^*(E))\cong\ker(A_E). 
\end{equation} 

For an arbitrary graph $E$ it is now also possible to determine the 
ideal structure of $C^*(E)$ and, in particular, its primitive ideal 
space \cite{bhrs,hs3}. 

\section{Quantum lens spaces} 

In \cite{vs}, Vaksman and Soibelman defined and analyzed $q$-analogues 
of the odd-dimensional spheres, as homogeneous spaces of the 
quantum special unitary groups of Woronowicz \cite{w2}. 
For $q\in(0,1]$, the $C^*$-algebra $C(S_q^{2n-1})$ is 
therein identified with a universal $C^*$-algebra generated by 
$n$ elements $z_1,\ldots,z_n$, subject to certain relations. 
In the classical case $q=1$, these generators $\{z_j\}$ are  
the coordinate functions for $S^{2n-1}\subset\C^n$. Hence the elements 
of the $*$-algebra $\Pol(S^{2n-1}_q)$, generated algebraically by 
$z_1,\ldots,z_n$, play the role of polynomials on the quantum 
sphere $S^{2n-1}_q$. As $C^*$-algebras, $C(S^{2n-1}_q)$ are all 
isomorphic for $q\in(0,1)$. The Vaksman-Soibelman 
relations make sense for $q=0$ as well and lead to an isomorphic 
$C^*$-algebra. 

Imitating the classical construction one can define quantum 
analogues of (generalized) lens spaces as follows \cite{hs2}. Choose 
an integer $p\geq2$ and integers $m_1,\ldots,m_n$ relatively prime to $p$. 
Let $\theta=e^{2\pi i/p}$. There exists an order $p$ automorphism 
$\tilde{\Lambda}$ of $C(S^{2n-1}_q)$ determined by 
\begin{equation}\label{ltilde}
\tilde{\Lambda}(z_i):=\theta^{m_i}z_i. 
\end{equation} 
The $C^*$-algebra $C(L_q(p;m_1,\ldots,m_n))$ of continuous functions 
on the quantum lens space $L_q(p;m_1,\ldots,m_n)$ is, by definition, 
the fixed-point algebra for $\tilde{\Lambda}$:
\begin{equation}
C(L_q(p;m_1,\ldots,m_n)):=C(S^{2n-1}_q)^{\tilde{\Lambda}}. 
\end{equation} 
Thus, if $p=2$ (and necessarily $m_1=\cdots=m_n=1$ and $\theta=-1$) 
then the such defined quantum lens spaces coincide with the 
quantum real projective spaces studied in \cite{hs1}. 
The polynomial algebra on $L_q(p;m_1,\ldots,m_n)$ is defined as the 
fixed-point algebra of the restriction of $\tilde{\Lambda}$ to 
$\Pol(S^{2n-1}_q)$. 

\section{Relation with graph algebras}

The odd-dimensional quantum spheres 
of Vaksman and Soibelman correspond to graph $C^*$-algebras, as follows. 
Let $L_{2n-1}$ be the directed graph with vertices $\{v_1,\ldots,v_n\}$ 
and edges $\{e_{i,j}\;|\;i=1,\ldots,n,\;j=i,\ldots,n\}$, with 
$s(e_{i,j})=v_i$ and $r(e_{i,j})=v_j$. An explicit (on generators) 
isomorphism 
\begin{equation} 
\phi:C(S^{2n-1}_q)\longrightarrow C^*(L_{2n-1}) 
\end{equation} 
is given in \cite[Theorem 4.4]{hs1}. Thus $\Lambda:=\phi\tilde{\Lambda}
\phi^{-1}$ is an automorphism of $C^*(L_{2n-1})$ such that the fixed-point 
algebra $C^*(L_{2n-1})^\Lambda$ is isomorphic via $\phi$ to 
$C(L_q(p;m_1,\ldots,m_n))$. This isomorphism allows one to study 
our quantum lens spaces with the help of the extensively developed 
machinery of Cuntz-Krieger algebras. In particular, $\Lambda$ is a
quasi-free automorphism of the graph algebra $C^*(L_{2n-1})$ 
determined by 
\begin{equation}\label{l}
\Lambda(S_{e_{i,j}})=\theta^{m_i}S_{e_{i,j}}. 
\end{equation} 
It turns out (cf. \cite[Theorem 2.5]{hs2}) that $C^*(L_{2n-1})^\Lambda$ 
(and hence $C(L_q(p;m_1,\ldots,m_n))$) is itself isomorphic to a certain 
graph algebra. Furthermore, the same is true about the crossed product 
$C^*(L_{2n-1})\rtimes_\Lambda\Z_p$. Indeed, by virtue of the results of 
\cite{kp}, this crossed product is naturally $*$-isomorphic to the 
$C^*$-algebra of the skew-product graph $L_{2n-1}\times_c\Z_p$. Here 
$c$ is a $\Z_p$-valued labeling of the edges of $L_{2n-1}$ such that $c(e_{i,j})=m_i$. 
The corresponding skew-product graph has vertices $L_{2n-1}^0\times\Z_p$ and 
edges $L_{2n-1}^1\times\Z_p$, with $s(e_{i,j},m)=(v_i,m-m_i)$ and 
$r(e_{i,j},m)=(v_j,m)$. Since the action $\Lambda$ is saturated \cite{hs2},  
the fixed-point algebra and the crossed product are strongly Morita 
equivalent. Consequently, 
\begin{equation}\label{lens}
C(L_q(p;m_1,\ldots,m_n))\cong_M C^*(L_{2n-1}\times_c\Z_p). 
\end{equation}

\section{Principal actions} 

As expected, the $\Z_p$-action $\tilde{\Lambda}$ defining our 
quantum lens spaces has nice properties resembling those of its 
classical counterpart. In particular, it is principal in the sense of 
Ellwood \cite{e}. In the present article we only consider actions 
of compact groups, which are automatically proper. Thus an 
action of a compact group $\Gamma$ on a $C^*$-algebra ${\mathcal A}$ 
is principal in the sense of \cite[Definition 2.4]{e} if 
and only if the linear map  
\begin{eqnarray}\label{can}
\Phi:{\mathcal A}\otimes_{\rm alg}{\mathcal A} & \longrightarrow & 
{\mathcal A}\otimes C(\Gamma), \nonumber \\ 
\Phi:x\otimes y & \longmapsto & (x\otimes I)\delta(y), \label{phi}
\end{eqnarray}
has norm-dense range. Here $\delta:{\mathcal A}\longrightarrow 
{\mathcal A}\otimes C(\Gamma)$ is the corresponding coaction.  

\begin{prop}
The action $\tilde{\Lambda}$ of $\Z_p$ on $C(S^{2n-1}_q)$ defined by 
(\ref{ltilde}) is principal. 
\end{prop}
{\em Proof.} It sufficies to prove that the action $\Lambda$ of 
$\Z_p$ on $C^*(L_{2n-1})$ defined by (\ref{l}) is principal. 
The corresponding coaction $\delta:C^*(L_{2n-1})\longrightarrow 
C^*(L_{2n-1})\otimes C(\Z_p)$ is given on the generators by 
\begin{equation}
\delta(S_{e_{i,j}})=S_{e_{i,j}}\otimes\chi^{m_i},  
\end{equation}  
where $\chi$ is the character of $\Z_p$ such that $\chi(d)=\theta^d$. Since 
$\theta$ is a $p^{\rm th}$ primitive root of unity, this unitary element 
$\chi$ generates the algebra $C(\Z_p)$ of functions on $\Z_p$. 

For each $i=1,\ldots,n$ and $k=1,2,\ldots$ we have 
\begin{equation}\label{eqn1}
\Phi(S_{e_{i,i}}^{*k}\otimes S_{e_{i,i}}^k)=(S_{e_{i,i}}^{*k}\otimes I)
(S_{e_{i,i}}^k\otimes\chi^{km_i})=P_{v_i}\otimes\chi^{km_i}. 
\end{equation} 
Since $m_i$ is relatively prime to $p$, by assumption, (\ref{eqn1}) 
implies that the range of $\Phi$ contains $P_{v_i}\otimes C(\Z_p)$. 
Thus, it also contains $I\otimes C(\Z_p)$, since in $C^*(L_{2n-1})$ 
we have $I=\sum_i P_{v_i}$. Consequently, the map $\Phi$ is surjective. 
$\square$ 

\vspace{3mm}\noindent
If $m_1=\cdots=m_n=1$ then the action $\Lambda$ coincides with the restriction 
of the gauge action to 
the group of $p^{\rm th}$ roots of unity. Recall that for an arbitrary graph 
$E$ the gauge action $\gamma:\T\longrightarrow\Aut(C^*(E))$ is defined by 
\begin{equation}
\gamma_t(S_e)=tS_e, \;\;\; \gamma_t(P_v)=P_v, 
\end{equation}
for all $e\in E^1$, $v\in E^0$ and $t\in\T\subset\C$. The $\T$-action on 
$C(S^{2n-1}_q)=C^*(z_1,\ldots,z_n)$ such that 
$t\cdot z_i=tz_i$ yields the fixed-point algebra $C(\C P^{n-1}_q)$ \cite{vs}. 
Under the isomorphism $C(S^{2n-1}_q)\cong C^*(L_{2n-1})$ of \cite[Theorem 4.4]{hs1} 
this action is transported into the gauge action $\gamma$ on 
$C^*(L_{2n-1})$. This observation motivates our next proposition. 

\begin{prop}\label{gauge}
If $E$ is a directed graph such that each vertex emits finitely many edges, 
each vertex emits at least one edge and receives at least one, then 
the gauge action $\gamma:\T\longrightarrow\Aut(C^*(E))$ is principal. 
\end{prop}
{\em Proof.} Let $z$ denote the canonical generator of $C(\T)$. Then 
the coaction $\delta:C^*(E)\longrightarrow C^*(E)\otimes C(\T)$, 
corresponding to the gauge action $\gamma$, is defined by 
\begin{equation}\label{d}
\delta(S_e)=S_e\otimes z, \;\;\; \delta(P_v)=P_v\otimes I, 
\end{equation}
for all $e\in E^1$, $v\in E^0$. Since the powers of $z$ span a dense subspace of 
$C(\T)$ and finite sums of projections of the form $P_v$, $v\in E^0$, 
give rise to an approximate unit for $C^*(E)$,  
it sufficies to show that the image of the map $\Phi$ (defined in  
(\ref{phi})) contains $P_v\otimes z^k$ for each $v\in E^0$, $k\in\Z$. 
To this end, fix a vertex $v$ and a positive integer $k$. 
Since each vertex of $E$ receives at least one edge it follows that 
there exists a path $\alpha=\alpha_1\cdots\alpha_k$, with $\alpha_i\in E^1$, 
which ends in $v$. For $S_\alpha=S_{\alpha_1}\cdots S_{\alpha_k}$ we have  
\begin{equation}\label{kp}
\Phi(S_\alpha^*\otimes S_\alpha)=(S_\alpha^*\otimes I)\delta(S_\alpha)
=(S_\alpha^*\otimes I)(S_\alpha\otimes z^k)=P_v\otimes z^k. 
\end{equation} 
Since each vertex of $E$ emits finitely many edges and at least one, 
an inductive application of (\ref{vertex}) yields 
\begin{equation}\label{pv}
P_v=\sum_{s(\beta)=v,|\beta|=k}S_\beta S_\beta^*. 
\end{equation} 
The summation in (\ref{pv}) extends over all paths of length $k$ in $E$ 
which begin at $v$. Thus we have 
\begin{equation}\label{km}
\Phi\left(\sum_{s(\beta)=v,|\beta|=k}S_\beta\otimes S_\beta^*\right)=
\sum_{s(\beta)=v,|\beta|=k}(S_\beta\otimes I)(S_\beta^*
\otimes z^{-k})=P_v\otimes z^{-k}. 
\end{equation} 
Combining (\ref{d}), (\ref{kp}) and (\ref{km}) we conclude that 
$\Phi$ satisfies the required property. 
$\square$ 

\vspace{3mm}\noindent
On the other hand, one can show that if $E$ contains a vertex which 
does not emit any edges then the gauge action is not principal. 

\begin{coro}
The $\T$-action on $C(S^{2n-1}_q)=C^*(z_1,\ldots,z_n)$ such that 
$t\cdot z_i=tz_i$, whose fixed-point algebra is $C(\C P^{n-1}_q)$ 
\cite{vs}, is principal. 
\end{coro} 
{\em Proof.} Combine Proposition \ref{gauge} above with 
\cite[Theorem 4.4 and \S 4.3]{hs1}. 
$\square$ 

\vspace{3mm}
In the case $n=2$, the $C^*$-algebra $C(S^3_q)\cong C(SU_q(2))$ is isomorphic to 
$C^*(L_3)$, corresponding to the following graph (cf. \cite{w1,hs1}): 

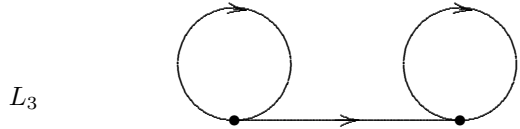
\[ \beginpicture
\setcoordinatesystem units <1.5cm,1.5cm>
\setplotarea x from -3 to 4, y from -0.1 to 1  
\put {$\bullet$} at 0 0
\put {$\bullet$} at 2 0 
\circulararc 360 degrees from 0 0 center at 0 0.5 
\circulararc 360 degrees from 2 0 center at 2 0.5
\setlinear 
\plot 0 0  2 0 / 
\put {$L_3$} [l] at -2 0.2
\arrow <0.235cm> [0.2,0.6] from 0 0.99 to 0.1 0.975 
\arrow <0.235cm> [0.2,0.6] from 2 0.99 to 2.1 0.975
\arrow <0.235cm> [0.2,0.6] from 0.9 0 to 1.1 0
\endpicture \]
The gauge action $\gamma$ on $C^*(L_3)$ is principal, with the fixed-point 
algebra isomorphic to the minimal unitization of the compacts. Differently 
interpreted, this setting corresponds on one hand to the quantum bundle 
$\T\rightarrow S^3_q\rightarrow \C P^1_q$ of \cite{vs}, 
and on the other hand to $\T\rightarrow SU_q(2)\rightarrow S^2_{q0}$ 
of \cite{p}. 

\vspace{5mm}\noindent
{\bf Acknowledgements.} I would like to thank Piotr M. Hajac for several 
very useful conversations on principal actions and other topics. I am grateful 
to the European Commission for partial support of my travel to Warsaw 
for this school/conference. It is a pleasure to thank Max-Planck-Institut 
f\"ur Mathematik (Bonn) and the Research Grants Committee (Newcastle) 
for their financial support. I would also like to thank 
Mathematisches Forschungsinstitut Oberwolfach, where this note was completed 
during my stay under the Research-in-Pairs programme.


\thebibliography{99} 

\bibitem{bhrs} T. Bates, J. H. Hong, I. Raeburn and W. Szyma\'nski, 
{\em The ideal structure of the $C^*$-algebras of infinite graphs}, 
Illinois J. Math., to appear. 

\bibitem{bprs} T. Bates, D. Pask, I. Raeburn and W. Szyma\'nski, 
{\em The $C^*$-algebras of row-finite graphs}, New York J. Math. 
{\bf 6} (2000), 307--324. 

\bibitem{c} J. Cuntz, {\em A class of $C^*$-algebras and topological 
Markov chains II: Reducible chains and the $Ext$-functor for $C^*$-algebras}, 
Invent. Math. {\bf 63} (1981), 25--40. 

\bibitem{ck} J. Cuntz and W. Krieger, {\em A class of $C^*$-algebras and 
topological Markov chains}, Invent. Math. {\bf 56} (1980), 251--268. 

\bibitem{e} D. A. Ellwood, {\em A new characterisation of principal 
actions}, J. Funct. Anal. {\bf 173} (2000), 49--60. 

\bibitem{hl} E. Hawkins and G. Landi, {\em Fredholm modules for quantum 
Euclidean spheres}, preprint, 2002. 

\bibitem{hs1} J. H. Hong and W. Szyma\'nski, {\em Quantum spheres and 
projective spaces as graph algebras}, Commun. Math. Phys., to appear.  

\bibitem{hs2} J. H. Hong and W. Szyma\'nski, {\em Quantum lens spaces and 
graph algebras}, preprint, 2001. 

\bibitem{hs3} J. H. Hong and W. Szyma\'nski, {\em The primitive ideal space 
of the $C^*$-algebras of infinite graphs}, preprint, 2002. 

\bibitem{kp} A. Kumjian and D. Pask, {\em $C^*$-algebras of directed graphs 
and group actions}, Ergodic Theory \& Dynamical Systems {\bf 19} 
(1999), 1503--1519. 

\bibitem{kprr} A. Kumjian, D. Pask, I. Raeburn and J. Renault, {\em 
Cuntz-Krieger algebras of directed graphs}, J. Funct. Anal. {\bf 144} 
(1997), 505--541. 

\bibitem{mt} K. Matsumoto and J. Tomiyama, {\em Noncommutative lens spaces}, 
J. Math. Soc. Japan {\bf 44} (1992), 13--41. 

\bibitem{p} P. Podle\'s, {\em Quantum spheres}, Lett. Math. Phys. 
{\bf 14} (1987), 193--202. 

\bibitem{rs} I. Raeburn and W. Szyma\'nski, {\em Cuntz-Krieger algebras of 
infinite graphs and matrices}, preprint, 1999. 

\bibitem{rtf} N. Yu. Reshetikhin, L. A. Takhtadzhyan, L. D. Faddeev, {\em Quantization 
of Lie groups and Lie algebras}, (Russian) Algebra i Analiz {\bf 1} (1989),
178--206; translation in Leningrad Math. J. {\bf 1} (1990), 193--225. 

\bibitem{vs} L. L. Vaksman and Y. S. Soibelman, {\em Algebra of functions 
on quantum $SU(n+1)$ group and odd dimensional quantum spheres}, 
Algebra-i-Analiz {\bf 2} (1990), 101--120. 

\bibitem{w1} S. L. Woronowicz, {\em Twisted ${\rm SU}(2)$ group. An example of 
a noncommutative differential calculus}, Publ. Res. Inst. Math. Sci. 
{\bf 23} (1987), 117--181. 

\bibitem{w2} S. L. Woronowicz, {\em Tannaka-Krein duality for compact 
matrix pseudogroups. Twisted $SU(N)$ groups}, Invent. Math. 
{\bf 93} (1988), 35--76.
\end{document}